\documentclass[10pt, conference, letterpaper]{IEEEtran}

\usepackage{amsmath}
\usepackage{amssymb}
\usepackage{amsbsy} % for bolding lower case greek letters
\usepackage{amsthm} % For theorems, lemma, definitions, propositions, conjectures, and proofs
\usepackage{bbm}
\usepackage{url}
\usepackage{array}
\usepackage{color}
\usepackage{multirow}
\usepackage{epsfig}
\usepackage{colortbl}
\usepackage [table]{xcolor}
\usepackage{graphicx}
\usepackage{float}
\usepackage{caption}
\usepackage{algorithm}
\usepackage{algpseudocode} % To use \Procedure command
\usepackage{mathtools} %for \DeclareMathOperator
\usepackage{siunitx} %for Angstroms 
\usepackage{subfigure}
\usepackage{booktabs} %for different tabular rules (lines)
\DeclareMathOperator*{\argmin}{arg\,min}
\DeclareMathOperator*{\minimize}{minimize}

\floatname{algorithm}{Algorithm}

\IEEEoverridecommandlockouts

\title{Spectrally Grouped Total Variation Reconstruction for Scatter Imaging Using ADMM}

\author{\IEEEauthorblockN{Ikenna~Odinaka\IEEEauthorrefmark{1}, Yan~Kaganovsky\IEEEauthorrefmark{1}, Joel~A.~Greenberg\IEEEauthorrefmark{1}, Mehadi~Hassan\IEEEauthorrefmark{1}, David~G.~Politte\IEEEauthorrefmark{2}, \\ Joseph A.~O'Sullivan\IEEEauthorrefmark{3}, Lawrence~Carin\IEEEauthorrefmark{1} and David J. Brady\IEEEauthorrefmark{1}}
\vspace{1.6mm}
\IEEEauthorblockA{\IEEEauthorrefmark{1} Department of Electrical and Computer Engineering, Duke University, Durham, NC 27708, USA.}
\IEEEauthorblockA{\IEEEauthorrefmark{2} Mallinckrodt Institute of Radiology, Washington University in Saint Louis, Saint Louis, MO 63110, USA.} 
\IEEEauthorblockA{\IEEEauthorrefmark{3} Department of Electrical and Systems Engineering, Washington University in Saint Louis, Saint Louis, MO 63130, USA.} 

\thanks{This paper appears in the proceedings of the 2015 IEEE Nuclear Science Symposium and Medical Imaging Conference (NSS/MIC).}% <-this % stops a space
}

\begin{document}

% make the title area
\maketitle

%% abstract
\begin{abstract}
We consider X-ray coherent scatter imaging, where the goal is to reconstruct momentum transfer profiles (spectral distributions) at each spatial location from multiplexed measurements of scatter. Each material is characterized by a unique momentum transfer profile (MTP) which can be used to discriminate between different materials. We propose an iterative image reconstruction algorithm based on a Poisson noise model that can account for photon-limited measurements as well as various second order statistics of the data. To improve image quality, previous approaches use edge-preserving regularizers to promote piecewise constancy of the image in the spatial domain while treating each spectral bin separately. Instead, we propose spectrally grouped regularization that promotes piecewise constant images along the spatial directions but also ensures that the MTPs of neighboring spatial bins are similar, if they contain the same material. We demonstrate that this group regularization results in improvement of both spectral and spatial image quality. We pursue an optimization transfer approach where convex decompositions are used to lift the problem such that all hyper-voxels can be updated in parallel and in closed-form.  The group penalty introduces a challenge since it is not directly amendable to these decompositions. We use the alternating directions method of multipliers (ADMM) to replace the original problem with an equivalent sequence of sub-problems that are amendable to convex decompositions, leading to a highly parallel algorithm. We demonstrate the performance on real data.
\end{abstract}

\section{Introduction \label{intro}}
We consider photon-limited X-ray coherent scatter imaging where the goal is to reconstruct momentum transfer profiles (spectral distributions) at each spatial location from multiplexed measurements of coherent X-ray scatter at a detector plane. Each material is characterized by a unique momentum transfer profile (MTP), so the reconstructed hyperspectral image can help discriminate between different materials at each spatial location.  Scatter imaging is ill-posed since the individual spectral components are not measured directly and regularization is desired in order to improve image quality. 
 
Standard approaches to reconstruction used in medical imaging often involve edge-preserving penalties to promote reconstructed images that are close to piecewise constant. One approach is to treat all dimensions as equal which implies a piecewise constant hyperspectral image in the spatio-spectral space. However this approach is not useful for coherent scatter imaging since MTPs are often spread out along the spectral axis. Another approach is to treat each spectral bin separately and promote images that are close to piecewise constant only in the spatial dimensions \cite{Ramani2011,Sawatzky2014,Odinaka2014}. This approach promotes uniformity across neighboring spatial locations (two neighboring spatial locations containing the same material should have the same MTP). However, this approach allows the spectral components of the image to change freely across spectral bins and may not account for possible correlations and smoothness across spectral components. 

To account for these situations, we study a type of regularization that groups the MTPs of spatially neighboring locations. It not only promotes piecewise constant images along the spatial directions, but also smoothness of the MTPs across spectral bins. Another motivation for using the group penalty stems from the fact that spectral resolution is better for higher spectral bins and most information for discriminating between different materials can be found in the lower spectral bins \cite{Kidane1999}. By grouping the spectral bins, we can effectively improve the resolution of the lower spectral bins, at the expense of the higher (non-discriminatory) ones. Variants of this form of regularization have appeared in previous works dealing with multichannel (including color) image denoising and recovery \cite{Chan01, Yang09}.

We propose an iterative image reconstruction algorithm based on a Poisson noise model which allows us to account for photon-limited measurements as well as various 2nd order statistics of the data, and the proposed group-TV penalty. 

To solve large problems, we pursue an optimization transfer approach where the original problem is replaced by an equivalent sequence of sub-problems that are simpler to optimize.  
Incorporating the proposed group penalty into this framework introduces a challenge since convex decompositions similar to De Pierro's trick \cite{DePierro1995} cannot be applied directly. To resolve this difficulty, we propose a new algorithm based on the alternating direction method of multipliers (ADMM) framework which breaks the problem into several simple sub-problems which can be solved using the mentioned optimization transfer approach, leading to a highly parallel algorithm. 

In summary, the main contributions of this paper are:
\begin{itemize}
\item We propose spectrally grouped TV regularization for hyperspectral image recovery in Poisson noise.
\item We develop a highly parallel algorithm based on ADMM and separable convex decompositions.
\item We study the performance of the new algorithm for coherent scatter imaging.
\end{itemize}

\section{The Proposed Model \label{problem_desc}}
The measurements $y_i$ from the X-ray coherent scatter system are modeled as independent Poisson random variables

\begin{equation}
\label{eq:poisson_model} 
y_i \sim \mbox{Poisson}(\sum_{j=1}^J a_{i,j} f_j + r_i), \qquad  i = 1, \dots, I 
\end{equation}
\noindent{where $I$ is the number of measurements, $J$ is the number of image voxels. $A \in \mathbb{R}_+^{I\times J}$ is the system matrix (forward model) with $a_{i,j}$ denoting the $ij$th entry. The column vector $\mathbf{f} \in \mathbb{R}_+^J$ is the lexicographical ordering of the hyperspectral image with $f_j$ denoting the $j$th entry. $\mathbf{r} \in \mathbb{R}_+^I$ are the background measurements which are assumed to be known with the $i$th entry denoted by $r_i$.} Let $J=S\times Q$, where $S$ is the number of spatial bins in the image and $Q$ is the number of spectral bins.

We propose an isotropic group TV penalty of the form
\begin{equation}
\label{eq:group_TV}
R(\mathbf{f}) = \sum_{s = 1}^{S}\sqrt{\sum_{q=1}^{Q} \sum_{p=1}^{|\mathcal{N}_s|} \left[ w_{s,q,p}\left(G_p \mathbf{f}\right)_{s,q}\right]^2},
\end{equation}
\noindent{where $G_p$ is a finite differencing matrix in the direction of the $p$th neighbor, $|\mathcal{N}_s|$ is the number of neighbors at the $s$th spatial bin, and $w_{s, q, p}$ is a weight to compensate for different physical units of the spectral and spatial dimensions, and different voxel sizes in each dimension. This penalty generalizes the multichannel TV penalty in \cite{Yang09} by allowing different weights for different neighboring directions $p$.}
We consider penalized log-likelihood estimation corresponding to \eqref{eq:poisson_model} and the penalty in \eqref{eq:group_TV}.  
 
Note that in the isotropic standard TV regularizer the sum over $q$ sits outside of the square root in \eqref{eq:group_TV}, which does not group the spectral bins and treats them independently.
In addition, note that \eqref{eq:group_TV} is not in a form that permits a direct application of convex decomposition methods. This motivates us to reformulate the problem and use ADMM to obtain a sequence of sub-problems that are amendable to convex decomposition.

\section{The Proposed Algorithm \label{algo_dev}}
We reformulate the penalized likelihood estimation as
\begin{equation}
\label{op:equiv_problem}
\begin{aligned} 
& \minimize_{\mathbf{f}, \mathbf{d}} & & L(\mathbf{f}) + \beta R(\mathbf{d}) \\
& \text{subject to} & & \mathbf{f} \geq \mathbf{0} \\
&\, & & \mathbf{d}_p = \mathbf{w}_p\odot G_p\mathbf{f}, \quad p=1,...,N
\end{aligned}
\end{equation}
\noindent {where $L$ is the log-likelihood corresponding to \eqref{eq:poisson_model}, $R$ is given in \eqref{eq:group_TV}, $\odot$ denotes element-wise multiplication}, and $\mathbf{w}_p$ and $\mathbf{d}_p$ are concatenations of $w_{s,q,p}$ and $d_{s,q,p}$ over all $(s, q)$. We assume for convenience that each voxel has the same number of spatial neighbors, i.e., $|\mathcal{N}_s|=N$ $\forall s$ (voxels at spatial boundaries are assumed to have additional neighbors outside of the domain with zero values) which is equivalent to Dirichlet boundary conditions. This requirement is not necessary and made here only to simplify notation. 

Using the split Bregman algorithm \cite{Goldstein2009} (which is equivalent to ADMM) to solve \eqref{op:equiv_problem} we have 
\begin{align}
&\mathbf{f}^{k + 1} = \argmin_{\mathbf{f}} L(\mathbf{f}) + \frac{\lambda}{2}\|\mathbf{d}^k - \mathbf{w}\odot\left(G\mathbf{f}\right) - \mathbf{c}^k\|_2^2 \label{eq:image_update} \\
&\mathbf{d}^{k + 1}  = \argmin_{\mathbf{d}} \beta R(\mathbf{d}) + \frac{\lambda}{2}\|\mathbf{d} - \mathbf{w}\odot\left(G\mathbf{f}^{k + 1}\right) - \mathbf{c}^k\|_2^2 \label{eq:splitting_variable_update} \\
&\mathbf{c}^{k + 1}  =  \mathbf{c}^{k} + \mathbf{w}\odot\left(G\mathbf{f}^{k + 1}\right) - \mathbf{d}^{k + 1} \label{eq:dual_update},
\end{align}
\noindent where $\mathbf{d}, \mathbf{w} \in \mathbb{R}^{J \times N}$, and $G$ are concatenations of all $\mathbf{d}_p$, $\mathbf{w}_p$, and $G_p$, for all $p$, respectively.
$\lambda > 0$ is the penalty parameter, and $\mathbf{c}$ is the set of scaled dual variables \cite{Boyd2011} or Bregman constants \cite{Goldstein2009}. The fully-parallel image update equations are obtained by using the fully-separable EM surrogate \cite{DePierro1995} for $L(\mathbf{f})$ and De Pierro's convexity trick \cite{DePierro1995} for the quadratic penalty in \eqref{eq:image_update}.

\section{Results \label{results}}
We consider an application of the proposed regularization scheme to a compressive fan beam system with energy-sensitive detectors, which has previously been described and characterized by Greenberg \textit{et al.} \cite{Greenberg2013snapshot, Greenberg2015}. The system consists of a 125 kvP X-ray source which has been collimated down to a small fan. The forward propagating X-rays illuminate a slice of an object volume, where scattering occurs in every direction. We focus on low-angle forward scatter, which is mainly due to coherent scattering, and collect scattered photons that are incident on a single energy-sensitive MULTIX \cite{Multix} linear detector array. The array has 128 pixels and 64 energy channels. The rays leaving the object volume are incident on a coded aperture \cite{Brady2013coding, Greenberg2015}, which allows for multiplexed measurement without the need for a rotating gantry. 

A 4 mm by 2 mm piece of Teflon was imaged by the system. The mean number of detected photons was approximately 26. For image recovery, we consider a 205 mm by 13.5 mm region within the fan plane, along the direction of ($z$) and perpendicular ($y$) to the central ray, respectively. The sampling interval in the $z$ and $y$ directions are 5 mm and 1.5 mm respectively. The 54 momentum transfer bins are uniformly sampled from 0.05 to 0.4475 inverse Angstroms (1/\si{\angstrom}). To permit a fair comparison, the regularization parameter $\beta$ was chosen for each regularizer, by sweeping over a range of values, and selecting the image with the most accurate spatial distribution; the spatial distribution is obtained by summing over the spectral bins. 

\begin{figure}
	\centering
	\subfigure[Group TV]{%
		\includegraphics[width=0.45\textwidth]{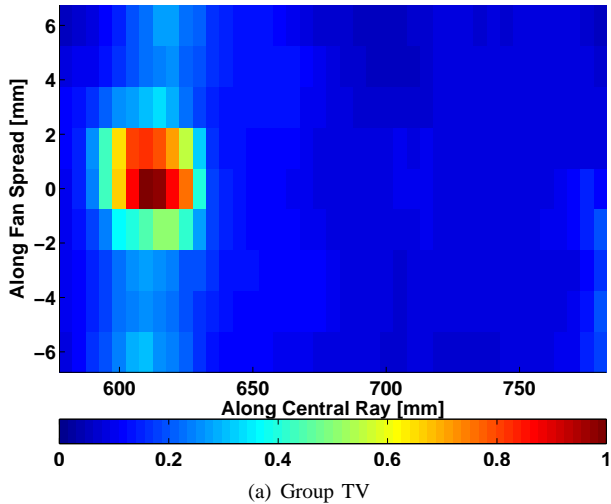}%
		\label{fig:recon_spatial_dist_teflon_group_TV}}%
	\\
	\subfigure[Standard TV]{%
		\includegraphics[width=0.45\textwidth]{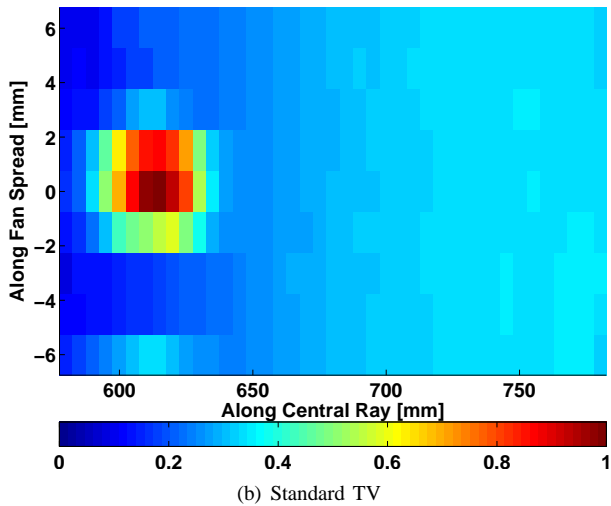}%
		\label{fig:recon_spatial_dist_teflon_standard_TV}}%    
	\caption{Reconstructed spatial distribution. Each spatial distribution was normalized to have a maximum value of 1. It was then non-linearly transformed to emphasize values close to 0 by taking the square root.}%
	\label{fig:recon_spatial_dist_teflon}%
\end{figure}

\begin{figure}
\centering
\includegraphics[width=0.45\textwidth]{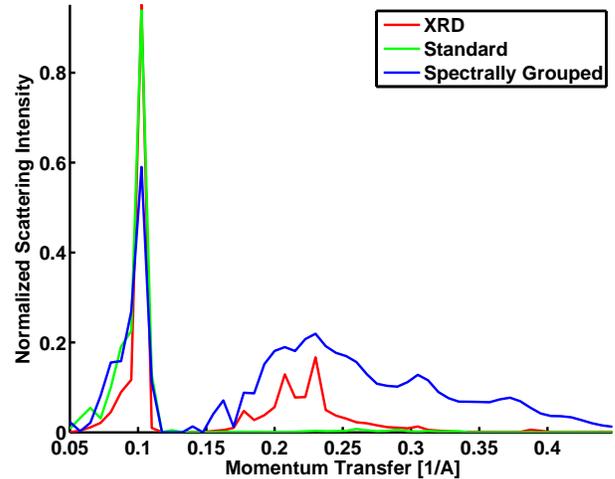}%
\caption{MTPs at the spatial location corresponding to the maximum of the spatial distribution for standard and group TV compared against a reference MTP obtained using an X-ray diffractometer (XRD). The MTPs are normalized to have unit energy.}%
\label{fig:mtps_teflon_xrd_standard_group}%
\end{figure}

Figure~\ref{fig:recon_spatial_dist_teflon} shows the spatial distributions of the reconstructed object using group and standard TV. Each spatial distribution was normalized to have a maximum value of 1. It was then non-linearly transformed to emphasize values close to 0 by taking the square root. From the figures, group TV appears to do a better job at keeping the spatial smoothing closer to where Teflon is located. Moreover, in Figure~\ref{fig:mtps_teflon_xrd_standard_group}, we see that group TV locates more of the characteristic peaks for Teflon than standard TV. However, standard TV more accurately recovers the first prominent peak. The reference MTP in Figure~\ref{fig:mtps_teflon_xrd_standard_group}, denoted as XRD, was obtained using an X-ray diffractometer \cite{Bruker}. As expected from spectral consensus, the reconstruction effort is distributed across the spectral bins for group TV, while standard TV dedicates most of its efforts around the prominent peak. 

\section{Conclusions \label{conc}}
We proposed a new spectrally grouped TV regularizer for recovering hyperspectral images in Poisson noise. We developed a highly parallel algorithm for image recovery based on ADMM and separable convex decompositions. In the experiment performed in X-ray coherent scatter imaging, the proposed regularizer outperforms the standard regularizer in both spectral and spatial image quality.

\bibliography{./References}
\bibliographystyle{IEEEtran}
\end{document}